\documentclass[10pt]{article}
\usepackage{amssymb, amsmath, url, graphicx, setspace, geometry}
\usepackage{calrsfs}
\usepackage{wasysym}

\def\3{\subset }
\def\4{\subseteq }
\def\<{\left<}
\def\>{\right>}

\def\bit{\begin{itemize}}
\def\eit{\end{itemize}}
\def\3{\subset }
\def\4{\subseteq }

\def\0{\leqno}

\def\barr{\begin{array}}
\def\earr{\end{array}}

\def\Z{{\rlap{$\kern2pt{\rm Z}$}{\rm Z}\,}}
\def\bld#1#2{{\buildrel{#1}\over{#2}}}
\def\st#1#2{{\mathrel{\mathop{#2}\limits_{#1}}{}\!}}
\def\stb#1#2#3{{\st{{#1}}{\bld{{#2}}{#3}}{}\!}}
\def\xmare#1#2{\stb{#1}{#2}{\mbox{\Large$\times$}}}

\title{\bf Some density results involving the average order of a finite group}
\author{Mihai-Silviu Lazorec}
\date{March 4, 2023}

\begin{document}

\maketitle

\begin{abstract}
Let $o(G)$ be the average of the element orders of a finite group $G$. A research topic concerning this quantity is understanding the relation between $o(G)$ and $o(H)$, where $H$ is a subgroup of $G$. Let $\mathcal{N}$ be the class of finite nilpotent groups and let $L(G)$ be the subgroup lattice of $G$. In this paper, we show that the set $\lbrace \frac{o(G)}{o(H)} \ | \ G\in\mathcal{N}, H\in L(G)\rbrace$ is dense in $[0, \infty)$. Other density results are outlined throughout the paper.
\end{abstract}



\section{Introduction}

Let $G$ be a finite group. In \cite{5}, A. Jaikin-Zapirain finds a super-logarithmic lower bound for the number of conjugacy classes $k(G)$ of $G$, when $G$ is nilpotent. More exactly, Theorem 1.1 of the same paper states that 
$$k(G)>10^{-4}\cdot\frac{\log_2\log_2n}{\log_2\log_2\log_2n}\cdot \log_2n,$$
where $G$ is a nilpotent group of order $n\geq 5$. One of the tools which plays a significant role in the proof of the above result is the so-called average order of $G$, i.e. the quantity
$$o(G)=\frac{1}{|G|}\sum\limits_{x\in G}|x|,$$
where $|x|$ denotes the order of an element $x\in G$. Among others, the author proves that $o(G)\geq o(Z(G))$, for any finite group $G$, and suggests that it would be interesting to further investigate the relation between the average order of $G$ and the average orders of its subgroups by answering the following question:\\

\textbf{Question 1.1.} \textit{Let $G$ be a finite ($p$-)group and let $N$ be a normal (abelian) subgroup of $G$. Is it true that $o(G)\geq o(N)^{\frac{1}{2}}?$}\\ 

Question 1.1 remained unanswered for nearly a decade. During 2021, E.I. Khukhro, A. Moret\' o and M. Zarrin published the paper \cite{6} which provides a negative answer to a generalized version of Jaikin-Zapirain's question. More exactly, Theorem 1.2 of \cite{6} states that given a real number $c>0$ and a prime number $p\geq \frac{3}{c}$, one can construct a $p$-group $G$ with a normal abelian subgroup $N$ such that $o(G)<o(N)^c$. Hence, for $c=\frac{1}{2}$, it is clear that there are counterexamples to Question 1.1. By following the notations in \cite{6}, these counterexamples are constructed by taking $G$ to be a semidirect product of a homocyclic group $U_s$ of exponent $p^{s}$, where $s=p+1$, and a so-called secretive $p$-group $P$ (see \cite{11} and Lemma 4.1 of \cite{6}), while $N$ is set to be $U_{s}$.  

Let $\mathcal{F}$ be the class of all finite groups, let $\mathcal{N}$ be the class of all finite nilpotent groups and let $L(G)$ be the subgroup lattice of a finite group $G$. For a subset $A$ of $\mathbb{R}$, we denote by $\overline{A}$ the closure of $A$ with respect to the usual topology $\tau_{\mathbb{R}}$ of $\mathbb{R}$. If we work with a different topology, say $\tau$, we denote the closure of $A$, with respect to $\tau$, by $\overline{A}_{\tau}$.
 
This paper also aims to investigate the relation between $o(G)$ and $o(H)$, where $H\in L(G)$, by studying the density of the set 
$$O_{\mathcal{C}}=\bigg\lbrace \frac{o(G)}{o(H)} \ \bigg| \ G\in\mathcal{C}, H\in L(G)\bigg\rbrace$$
in $[0, \infty)$, where $\mathcal{C}$ is a specific class of finite groups. We manage to show that $O_{\mathcal{F}}$ is dense in $[0, \infty)$ as a consequence of our main result which is even stronger and  states that:\\

\textbf{Theorem 1.2.} \textit{The set $O_{\mathcal{N}}$ is dense in $[0, \infty)$.}\\

An immediate consequence of Theorem 1.2 is obtained as follows. Let $\mathcal{C}$ be a class of finite groups such that $\mathcal{N}\subseteq \mathcal{C}$. Then $O_{\mathcal{N}}\subseteq O_{\mathcal{C}}\subseteq [0, \infty)$, so $\overline{O_{\mathcal{N}}}\subseteq \overline{O_{\mathcal{C}}}\subseteq \overline{[0, \infty)}$. Since $[0, \infty)$ is a closed set and $\overline{O_{\mathcal{N}}}=[0, \infty)$, we get:\\

\textbf{Corollary 1.3.} \textit{Let $\mathcal{C}$ be a class of finite groups such that $\mathcal{N}\subseteq  \mathcal{C}$. Then $O_{\mathcal{C}}$ is dense in $[0, \infty)$. In particular, $O_{\mathcal{F}}$ is dense in $[0, \infty)$.}\\

We end the introduction by mentioning that the average order of a finite group $G$ may be also expressed as
$$o(G)=\frac{\psi(G)}{|G|},$$
where $\psi(G)=\sum\limits_{x\in G}|x|$ is the sum of element orders of $G$. During the last years, there was a growing interest in investigating this invariant. We refer the reader to \cite{4} for a recent survey including relevant results concerning the sum of element orders of a finite group. 
  
\section{Proof of Theorem 1.2 and other results concerning the density of some sets}

As it was suggested in the first section, to obtain the density of $O_{\mathcal{F}}$ in $[0, \infty)$, it would be sufficient to find a class of groups $\mathcal{C}\subseteq\mathcal{F}$ such that $\overline{O_{\mathcal{C}}}=[0, \infty)$. It is clear that $\overline{O_{\mathcal{C}}}\subseteq [0, \infty)$. So, once we choose a candidate for the class $\mathcal{C}$, it suffices to show that each $a\in [0, \infty)$ is an adherent point of $O_{\mathcal{C}}$, i.e. there is a sequence of groups $(G_n)_{n\geq 1}\subset \mathcal{C}$ and a corresponding sequence $(H_n)_{n\geq 1}$, where $H_n\in L(G_n)$, for all $n\geq 1$, such that $$\displaystyle\lim_{n\to\infty}\frac{o(G_n)}{o(H_n)}=a.$$ 

Our candidate for $\mathcal{C}$ is $\mathcal{N}$ and, in what follows, we justify this option. To expand our reasoning, we include the following preliminary result which is a consequence of the Proposition outlined on p. 863 of \cite{7}.\\

\textbf{Lemma 2.1.} \textit{Let $(x_n)_{n\geq 1}$ be a sequence of positive real numbers such that 
$$\displaystyle\lim_{n\to\infty}x_n=0 \text{ \ and \ } \sum\limits_{n=1}^{\infty}x_n=\infty.$$
Then the set containing the sums of all finite subsequences of $(x_n)_{n\geq 1}$ is dense in $[0, \infty)$.}\\

We denote the $n$th prime number by $p_n$. Lemma 2.1 is the main tool that is going to be used to show that each $a\in [1, \infty)$ is an adherent point of $O_{\mathcal{N}}$. Once this is done, it remains to cover the points $a\in [0, 1)$. For the first part, the main idea is to apply Lemma 2.1 for a sequence $(x_n)_{n\geq 1}$, where $x_n=\ln \frac{o(\widetilde{G_n})}{o(\widetilde{H_n})}$. We are going to show that some suitable candidates for $(\widetilde{G_n})_{n\geq 1}$ and $(\widetilde{H_n})_{n\geq 1}$, such that the sequence $(x_n)_{n\geq 1}$ defined above satisfies the hypotheses of Lemma 2.1, are $\widetilde{G_n}=C_{p_n}^m$ and $\widetilde{H_n}=C_{p_n}^{m-1}$ for a fixed integer $m\geq 2$ (see the proof of Claim 2.5 below). Consequently, by applying Lemma 2.1 and some calculus properties, we deduce that there exists a sequence $(G_n)_{n\geq 1}$ of finite abelian groups and a corresponding sequence of subgroups $(H_n)_{n\geq 1}$ such that
$$\displaystyle\lim_{n\to\infty}\frac{o(G_n)}{o(H_n)}=a\in [1, \infty).$$
This means that $[1, \infty)\subseteq \overline{O_{\mathcal{A}}}$, where $\mathcal{A}$ is the class of finite abelian groups. The reverse inclusion also holds because $\frac{o(G)}{o(H)}\geq 1$ for any finite abelian group $G$ and any $H\in L(G)$. Indeed, since $G$ is self dual (see Chapter 8 of \cite{9} or \cite{2}), we know that for any $H\in L(G)$, there is $K\in L(G)$ such that $H\cong \frac{G}{K}$. Hence,
$$\frac{o(G)}{o(H)}=\frac{o(G)}{o(\frac{G}{K})}=\frac{1}{|K|}\cdot \frac{\psi(G)}{\psi(\frac{G}{K})}=\frac{\sum\limits_{x\in G}|x|}{\sum\limits_{x\in G}|xK|}\geq 1,$$
so $O_{\mathcal{A}}\subseteq [1, \infty)$ and this leads to $\overline{O_{\mathcal{A}}}\subseteq [1, \infty)$. Thus, we state the following result.\\

\textbf{Corollary 2.2.} \textit{The set $O_{\mathcal{A}}$ is dense in $[1, \infty)$.}\\

We mention that Corollary 2.2 also holds if we replace $\mathcal{A}$ with a class $\widetilde{\mathcal{C}}$ of finite groups such that $O_{\mathcal{A}}\subseteq O_{\widetilde{\mathcal{C}}}\subseteq [1, \infty)$. 

Finally, concerning the adherence property of the points $a\in [0, 1)$, we will mainly work with sequences formed of specific direct products of finite $p$-groups. Each such direct product has two main components: one is abelian, while the other one is a counterexample to Question 1.1 (see the proof of Claim 2.7 below). All finite groups that were highlighted in the last paragraphs are nilpotent and this consequently explains why our choice for $\mathcal{C}$ is $\mathcal{N}$.\\

The following preliminary result includes some number theoretic and calculus properties which are going to be used further.\\

\textbf{Lemma 2.3.} \textit{\begin{itemize}
\item[i)] Let $G_1$ and $G_2$ be finite groups. If $(|G_1|, |G_2|)=1$, then $$o(G_1\times G_2)=o(G_1)\cdot o(G_2).$$
\item[ii)] $$\sum\limits_{n=1}^{\infty}\frac{1}{p_n}=\infty.$$
\item[iii)] Let $(x_n)_{n\geq 1}, (y_n)_{n\geq 1}$ be sequences of positive real numbers. If 
$$\displaystyle\lim_{n\to\infty}\frac{x_n}{y_n}\in (0, \infty),$$
then the series $\sum\limits_{n=1}^{\infty}x_n$ and $\sum\limits_{n=1}^{\infty}y_n$ have the same nature.
\item[iv)] Let $(X, \tau)$ and $(Y, \tau')$ be topological spaces, let $f:X\longrightarrow Y$ be a continuous function and let $A, B\subseteq X$. If $\overline{A}_{\tau}=\overline{B}_{\tau}$, then $\overline{f(A)}_{\tau'}=\overline{f(B)}_{\tau'}$.
\end{itemize}}

Concerning the previous lemma, we mention that item \textit{i)} states that the average order is a multiplicative function. This is a consequence of the multiplicativity of the sum of element orders (see Lemma 2.1 of \cite{1}). A short proof of item \textit{ii)} may be found in \cite{8}. For item \textit{iii)}, one can check Theorem 10.9 of \cite{3}, while item \textit{iv)} is easily obtained using the characterization of the continuity of a function in terms of closure (see Proposition 6.12 of \cite{10}).

Let $I=[1, \infty)$. Denote by $\tau_I$ the subspace topology on $I$. For a subset $A$ of $\mathbb{R}$, the closure of $A$ with respect to $\tau_I$ is $\overline{A}_{\tau_I}=\overline{A}\cap I.$ By Corollary 2.2, we have $\overline{O_{\mathcal{A}}}=I$. We deduce that
\begin{equation}\label{r8}
\overline{O_{\mathcal{A}}}_{\tau_I}=\overline{I}_{\tau_I}.
\end{equation}
Since the function
$$f:(I, \tau_I)\longrightarrow (\mathbb{R}, \tau_{\mathbb{R}}), \text{ \ given by \ } f(x)=\frac{1}{x}, \ \forall \ x\in I,$$
is continuous, by Lemma 2.3, \textit{iv)}, and (\ref{r8}), we get
$$\overline{\bigg\lbrace \frac{o(H)}{o(G)} \ \bigg | \ G\in \mathcal{A}, H\in L(G)\bigg\rbrace}=\overline{(0, 1]}=[0, 1].$$
Therefore, one can state the following result.
\\

\textbf{Corollary 2.4.} \textit{The set
$$\bigg\lbrace \frac{o(H)}{o(G)} \ \bigg | \ G\in \mathcal{A}, H\in L(G)\bigg\rbrace$$
is dense in $[0, 1]$.}\\

We proceed now with the proof of the main result.\\

\textbf{Proof of Theorem 1.2.} Recall that $p_n$ denotes the $n$th prime number. We are going to complete some preliminary steps towards achieving our goal.\\

\textbf{Claim 2.5.} \textit{Let $m\geq 2$ be an integer. The set
$$\Bigg\lbrace \frac{o\big(\xmare{n\in I}{ }C_{p_n}^m\big)}{o\big(\xmare{n\in I}{ }C_{p_n}^{m-1}\big)} \ \Bigg| \ I\subset\mathbb{N}^*, |I|<\infty \Bigg\rbrace$$
is dense in $[1, \infty)$.}\\

\textbf{Proof.} Consider the sequence $(x_n)_{n\geq 1}$, where $x_n=\ln\frac{o(C_{p_n}^m)}{ o(C_{p_n}^{m-1}) }$, for all $n\geq 1$. We have
$$x_n=\ln\frac{p_n^{m+1}-p_n+1}{p_n^{m+1}-p_n^m+p_n^{m-1}}\in (0, \infty).$$
As $n$ approaches infinity, we get
\begin{equation}\label{r1}
\displaystyle\lim_{n\to\infty}x_n=\ln 1=0.
\end{equation}
Further, take the sequence $(y_n)_{n\geq 1}$ given by $y_n=\frac{1}{p_n}$, for all $n\geq 1$. Then
$$\displaystyle\lim_{n\to\infty}\frac{x_n}{y_n}=\displaystyle\lim_{n\to\infty}\bigg(p_n\cdot\ln\frac{p_n^{m+1}-p_n+1}{p_n^{m+1}-p_n^m+p_n^{m-1}}\bigg)=1\in (0, \infty).$$
By Lemma 2.3, \textit{ii)}, \textit{iii)}, we have
\begin{equation}\label{r2}
\sum\limits_{n=1}^{\infty}x_n=\infty.
\end{equation}
According to (\ref{r1}) and (\ref{r2}), the sequence $(x_n)_{n\geq 1}$ satisfies the hypotheses of Lemma 2.1. Hence, we have
\begin{equation}\label{r3}
\overline{\bigg\lbrace\sum\limits_{n\in I}x_n \ \bigg| \ I\subset\mathbb{N}^*, |I|<\infty \bigg\rbrace}=[0, \infty)\Longleftrightarrow \overline{\bigg\lbrace\ln\bigg(\prod_{n\in I}\frac{o(C_{p_n}^m)}{o(C_{p_n}^{m-1})}\bigg) \ \bigg| \ I\subset\mathbb{N}^*, |I|<\infty \bigg\rbrace}=[0, \infty).
\end{equation}
Since, by Lemma 2.3, \textit{i)}, the average order is a multiplicative function, (\ref{r3}) becomes
\begin{equation}\label{r4}
\overline{\Bigg\lbrace\ln\frac{o\big(\xmare{n\in I}{ }C_{p_n}^m\big)}{o\big(\xmare{n\in I}{ }C_{p_n}^{m-1}\big)} \ \Bigg| \ I\subset\mathbb{N}^*, |I|<\infty \Bigg\rbrace}=[0, \infty).
\end{equation}
Finally, since $$exp:(\mathbb{R}, \tau_{\mathbb{R}})\longrightarrow (\mathbb{R}, \tau_{\mathbb{R}}), \text{ \ given by \ } exp(x)=e^x, \forall \ x\in\mathbb{R},$$
is continuous and (\ref{r4}) highlights the equality of two closed sets of $(\mathbb{R}, \tau_{\mathbb{R}})$, we apply Lemma 2.3, \textit{iv)}, to finish the proof of our claim, i.e.
$$\overline{\Bigg\lbrace\frac{o\big(\xmare{n\in I}{ }C_{p_n}^m\big)}{o\big(\xmare{n\in I}{ }C_{p_n}^{m-1}\big)} \ \Bigg| \ I\subset\mathbb{N}^*, |I|<\infty \Bigg\rbrace}=[1, \infty).$$
\hfill\rule{1,5mm}{1,5mm}\\

\textbf{Claim 2.6.} \textit{Let $m\geq 2$ and let $J$ be a finite non-empty subset of $\mathbb{N}^*$. The set
$$\Bigg\lbrace \frac{o\big(\xmare{n\in I}{ }C_{p_n}^m\big)}{o\big(\xmare{n\in I}{ }C_{p_n}^{m-1}\big)} \ \Bigg| \ I\subset\mathbb{N}^*\setminus J, |I|<\infty \Bigg\rbrace$$
is dense in $[1, \infty)$.}\\

\textbf{Proof.} This is obtained by repeating the proof of Claim 2.5 for the sequence $(\widetilde{x_n})_{n\in\mathbb{N}^*\setminus J}$, where $\widetilde{x_n}=\ln\frac{o(C_{p_n}^m)}{ o(C_{p_n}^{m-1}) }$, for all $n\in \mathbb{N}^*\setminus J$. The same reasoning can be repeated since $(\widetilde{x_n})_{n\in\mathbb{N}^*\setminus J}$ is obtained by removing a finite number of terms from the original sequence $(x_n)_{n\geq 1}$ taken in the proof of Claim 2.5, so $(\widetilde{x_n})_{n\in\mathbb{N}^*\setminus J}$ also satisfies the hypotheses of Lemma 2.1.
\hfill\rule{1,5mm}{1,5mm}\\

\textbf{Claim 2.7.} \textit{Any $a\in [0, 1)$ is an adherent point of $O_{\mathcal{N}}$.}\\

\textbf{Proof.} Suppose that $a=0$. As we outlined in the first section, for $n\geq 4$ (i.e. for a prime greater than or equal to 7), if we take $G_n=U_{s_n}P_n$ to be a semidirect product of a homocyclic group $U_{s_n}$ of exponent $p_n^{s_n}$, where $s_n=p_n+1$, and a secretive $p_n$-group $P_n$, one has  $o(G)<o(U_{s_n})^{\frac{1}{2}}$. According to the proof of Theorem 1.2 of \cite{6}, the following inequalities hold:
$$o(G_n)< p_n^3 \text{ \ and \ } o(U_{s_n})\geq p_n^{p_n}, \ \forall \ n\geq 4.$$
Hence,
\begin{equation}\label{r5}
\frac{o(G_n)}{o(U_{s_n})}<\frac{p_n^3}{p_n^{p_n}}.
\end{equation}
As $n$ approaches infinity, (\ref{r5}) leads us to
\begin{equation}\label{r6}
\displaystyle\lim_{n\to\infty}\frac{o(G_n)}{o(U_{s_n})}=0,
\end{equation}
so $a=0$ is an adherent point of $O_{\mathcal{N}}$.

Let $a\in (0, 1)$. By (\ref{r6}), there is a sufficiently large $N$ such that $a\geq\frac{o(G_N)}{o(U_{s_N})}$. Consequently, $a\cdot \frac{o(U_{s_N})}{o(G_N)}\in [1, \infty)$. If we take $J=\lbrace N\rbrace$ in Claim 2.6, it follows that there is a sequence of finite abelian groups $(\widetilde{G_n})_{n\geq 1}$ and a corresponding sequence $(\widetilde{H_n})_{n\geq 1}$, where $\widetilde{H_n}\in L(\widetilde{G_n})$ for all $n\geq 1$, such that
\begin{equation}\label{r7}
\displaystyle\lim_{n\to\infty}\frac{o(\widetilde{G_n})}{o(\widetilde{H_n})}=a\cdot \frac{o(U_{s_N})}{o(G_N)}.
\end{equation}
Finally, we consider the sequences $(G_N\times \widetilde{G_n})_{n\geq 1}$ and $(U_{s_N}\times\widetilde{H_n})_{n\geq 1}$. Note that $(|G_N|, |\widetilde{G_n}|)=(|U_{s_N}|, |\widetilde{H_n}|)=1$, for all $n\geq 1$. Hence, by Lemma 2.3, \textit{i)}, and (\ref{r7}), we conclude that
$$\displaystyle\lim_{n\to\infty}\frac{o(G_N\times \widetilde{G_n})}{o(U_{s_N}\times\widetilde{H_n})}=\frac{o(G_N)}{o(U_{s_N})}\cdot\displaystyle\lim_{n\to\infty}\frac{o(\widetilde{G_n})}{o(\widetilde{H_n})}=a.$$
Hence, any $a\in (0, 1)$ is also an adherent point of $O_{\mathcal{N}}$ and this concludes the proof of our claim.
\hfill\rule{1,5mm}{1,5mm}\\

By Claims 2.5 and 2.7, it follows that $[0, \infty)\subseteq \overline{O_{\mathcal{N}}}$. Since the reverse inclusion also holds, the proof of Theorem 1.2 is complete.
\hfill\rule{1,5mm}{1,5mm}\\

We end our paper by posing a question concerning the class $\mathcal{P}$ of finite $p$-groups. If the answer would be affirmative, our main result would also follow since $\mathcal{P}\subset\mathcal{N}$.\\

\textbf{Question 2.8.} \textit{Is the set $O_{\mathcal{P}}$ dense in $[0, \infty)$?}

\bigskip\noindent {\bf Acknowledgements.} The author is grateful to the reviewers for their remarks which improve the previous version of the paper. This work was supported by a grant of the "Alexandru Ioan Cuza" University of Iasi, within the Research Grants program, Grant UAIC, code GI-UAIC-2021-01.

\vspace*{3ex}
\small
\hfill
\begin{minipage}[t]{7cm}
Mihai-Silviu Lazorec \\
Faculty of  Mathematics \\
"Al.I. Cuza" University \\
Ia\c si, Romania \\
e-mail: {\tt silviu.lazorec@uaic.ro}
\end{minipage}
\end{document}